\documentclass[12pt]{article}
\usepackage{amsfonts,amsmath,amsxtra}

\newtheorem{theorem}{Theorem}
\newtheorem{corollary}{Corollary}
\newtheorem{lemma}{Lemma}

\newenvironment{definition}
{\smallskip\noindent{\bf Definition\/}:}{\smallskip\par}
\newenvironment{proposition}
{\smallskip\noindent{\bf Proposition\/}.}{\smallskip\par}

\newenvironment{proof}
{\noindent{\bf Proof\/}.}{{ $\Box$}\smallskip\par}

\title{On the motivic measure on the space of functions.}
\author{E. Gorsky\footnote{Supported by the grants NSh-4719.2006.1, RFBR-04-01-00762.}}
\date{}

\begin{document}

\maketitle

\begin{abstract}

Motivic measure on the space of functions was introduced by
Campillo, Delgado and Gusein-Zade as an analog of
the motivic measure on the space of arcs .
In this paper we prove that the measure on the space of functions
can be related to the motivic measure on the space of arcs
by a factor, which can be defined explicitly in geometric terms.
This provides a possibility to rewrite motivic integrals over the
space of functions as integrals over the union 
of all symmetric powers of the space of arcs. 

\end{abstract}

\section{Introduction}

Motivic measure on the space of functions was introduced by
Campillo, Delgado and Gusein-Zade in \cite{al2} (see also \cite{alex}) as an analog of
the motivic measure on the space of arcs (see e.g. \cite{dl}).

Consider the map
$$Z:\mathbb{P}\mathcal{O}_{\mathbb{C}^2,0}\rightarrow \sqcup_{k}S^k\mathcal{B},$$
which maps a germ of a function $f$ on the plane (defined up to multiplication by a constant) to the collection of uniformisations 
of components of the germ 
of the curve 
$\{f=0\}$
up to automorphism of $(\mathbb{C},0)$ (here $\mathcal{B}$ denotes the quotient space of the space of arcs on $(\mathbb{C}^2,0)$  by the natural action of $Aut_{\mathbb{C},0}$).

In this paper we prove that the measure on the space of functions
can be related with the motivic measure on the space of arcs.
Namely, for every measurable set $M\subset \sqcup_{k}S^k\mathcal{B}$ one has
$$\mu(Z^{-1}(M))=\int_{M}\mathbb{L}^{-\delta(\gamma)-v(\gamma)}d\mu,$$  where $v(\gamma)$ is the order of the arc $\gamma$ at the origin, and $\delta(\gamma)$ is the number of self-intersection points of a generic deformation of $\gamma$.
This provides a possibility to rewrite motivic integrals over the
space of functions as integrals over the union 
of all symmetric powers of the space of arcs. 

The Euler characteristic morphism could not be applied directly 
to this equation. An analogous 
 correspondence formula for the Euler characteristics is proved
in section 5.  

\section{ Motivic measures}

Let $\mathcal{L}=\mathcal{L}_{\mathbb{C}^2,0}$ be the space
of arcs at the origin on the plane. 
It is the set of pairs $(x(t),y(t))$ of formal power series 
(without degree 0 term). 
Let $\mathcal{L}_n$ be the space of $n$-jets of such arcs, let $\pi_{n}:\mathcal{L}\rightarrow \mathcal{L}_{n}$ be the natural projection.

Let $K_0(Var_{\mathbb{C}})$ be the Grothendieck ring of quasiprojective
complex algebraic varieties. It is generated by the isomorphism classes
of complex quasiprojective algebraic varieties modulo the relations
 $[X]=[Y]+[X\setminus Y],$ where $Y$ is a Zariski closed subset of  $X$. 
Multiplication is given by the formula $[X]\cdot [Y]=[X\times Y].$
Let $\mathbb{L}\in K_0(Var_{\mathbb{C}})$ be the class of the complex affine line. 

The Euler characteristic provides
a ring homomorphism 
$$\chi: K_0(Var_{\mathbb{C}})\rightarrow \mathbb{Z}.$$

Let $S^{k}X$ denote the $k$-th symmetric power of the variety
$X$.

\begin{lemma}(\cite{powers})
For every $k$ and $m$ $[S^{k}\mathbb{C}^{m}]=\mathbb{L}^{km}.$
Moreover, for every $X$ $$[S^{k}(X\times\mathbb{C}^{m})]=\mathbb{L}^{km}\cdot S^{k}X.$$

\end{lemma}

Consider the ring $K_0(Var_{\mathbb{C}})[\mathbb{L}^{-1}]$
with the following filtration: $F_k$ is generated by the elements
of the type $[X]\cdot [\mathbb{L}^{-n}]$ with $n-\dim X\ge k$. Let $\mathcal{M}$
be the completion of the ring $K_0(Var_{\mathbb{C}})[\mathbb{L}^{-1}]$ corresponding
to this filtration.

On an algebra of subsets of $\mathcal{L}$ Kontsevich, and later Denef and Loeser (\cite{dl}) 
constructed a measure $\mu$
with values in the ring  $\mathcal{M}$. 

A subset $A\subset \mathcal{L}$ is said to be cylindric if there exist  $n$ and a constructible set $A_n\subset \mathcal{L}_{n}$ such that
$A=\pi_{n}^{-1}(A_n)$. For the cylindric set $A$ define
$$\mu(A)=[A_n]\cdot \mathbb{L}^{-2n}.$$
It was proved in \cite{dl}, that this measure can be extended to an additive measure on a suitable algebra of subsets in $\mathcal{L}$.    

A function $f:\mathcal{L}\rightarrow G$ with values in an abelian group $G$ is called simple, if its image is countable or finite, and for every $g\in G$ the set $f^{-1}(g)$ is measurable. Using this measure, one can define in the natural way the (motivic) integral for simple functions on $\mathcal{L}$ as
$\int_{\mathcal{L}}f d\mu=\sum_{g\in G} g\cdot\mu(f^{-1}(g)),$
if  the right hand side sum converges in  $G\otimes \mathcal{M}.$

Note that for cylindric sets the Euler characteristic can be well defined by the formula $\chi(A)=\chi(A_n)$. This gives a $\mathbb{Z}-$valued measure on the algebra of cylindric sets.
However, it cannot be extended to the algebra of measurable sets.  
This measure provides a notion of an integral with respect to the Euler characteristic for functions on $\mathcal{L}$ with cylindric level sets. It is clear that for such functions
$$\chi(\int_{\mathcal{L}}f d\mu)=\int_{\mathcal{L}}f d\chi .$$  

We will use simple functions
$v_{x}=Ord_{0}\,x(t), v_y=Ord_{0}\,y(t)$ and $v=\min\{v_x,v_y\}$,
defined for a curve $\gamma(t)=(x(t),y(t)).$

Campillo, Delgado and Gusein-Zade (\cite{alex}) 
constructed an analogous measure on the space
 $\mathcal{O}_{\mathbb{C}^2,0}=\mathcal{O}$ of germs of analytic functions on the plane at the origin.
Let $j_{n}(\mathcal{O})$ be the space of $n$-jets  of functions 
from $\mathcal{O}$.
A subset $A\subset \mathcal{O}$ is said to be cylindric if there exist  $n$ and a constructible set $A_n\subset j_{n}(\mathcal{O})$ such that
$A=j_{n}^{-1}(A_n)$. For the cylindric set $A$ define
$$\mu(A)=[A_n]\cdot \mathbb{L}^{1-{(N+1)(N+2)\over 2}}.$$ 
In the same way one can define the motivic integral over the space of functions.

Below we will need motivic measures
on the quotient spaces of the space of curves by the action of the groups
$\mathbb{C}^{*}$ and $Aut_{\mathbb{C},0}$ (the group of automorphisms of $(\mathbb{C},0)$) introduced in \cite{nonpar}.

On the space $\mathcal{L}$  there is
a natural action of $\mathbb{C}^{*}$, defined by
the formula $a\cdot\gamma(t)=\gamma(at),a\in\mathbb{C}^{*}.$
Let $\mathcal{L}^{*}=\mathcal{L}\setminus\{0\}.$

\begin{proposition}(\cite{nonpar})Let $p:\mathcal{L}^{*}\rightarrow\mathcal{L}^{*}/\mathbb{C}^{*}$ be
the factorization map. Then $\mu(p^{-1}(X))=(\mathbb{L}-1)\mu(X).$
\end{proposition} 

For an arc $\gamma:(\mathbb{C},0)\rightarrow(\mathbb{C}^2,0)$
and $h\in Aut_{\mathbb{C},0}$ let $h\cdot\gamma(t)=\gamma(h^{-1}(t)).$
It defines an action of the group $Aut_{\mathbb{C},0}$ on the space
of curves.

\begin{definition}(\cite{nonpar})
An orbit of this action is called a branch.
The space of branches will be denoted as $\mathcal{B}$.
\end{definition}

Let $p_0:\mathcal{L}\rightarrow\mathcal{B}$ be
the factorization map.

\begin{proposition}(\cite{nonpar})
Let $\psi$ be a simple integrable function on  $\mathcal{B}$.
Then $$\int_{\mathcal{L}}(p_0^{*}\psi)\cdot\mathbb{L}^{v}d\mu=(\mathbb{L}-1)\int_{\mathcal{B}}\psi d\mu.$$
\end{proposition}




\section{ Power structures}

The notion of the power structure over a (semi)ring was introduced
by S. Gusein-Zade, I. Luengo and A. Melle-Hernandez in \cite{powers}.

\begin{definition}
Power structure over a ring $R$
is a map $$(1+tR[[t]])\times R\rightarrow 1+tR[[t]]:(A(t),m)\mapsto (A(t))^{m},$$
satisfying the following properties:

1.$(A(t))^{0}=1,$

2.$(A(t))^{1}=A(t),$

3.$((A(t)\cdot B(t))^{m}=((A(t))^{m}\cdot((B(t))^{m},$

4.$(A(t))^{m+n}=(A(t))^{m}\cdot(A(t))^{n},$

5.$(A(t))^{mn}=((A(t))^{n})^{m},$

6.$(1+t)^{m}=1+mt+$ terms with higher degree,

7.$(A(t^{k}))^{m}=((A(t))^{m})|_{t\rightarrow t^{k}}.$
\end{definition}

For example, for $R=\mathbb{Z}$ the usual raising into a power  gives
a power structure, since if $A(t)\in 1+t\mathbb{Z}[t]$ and $m\in \mathbb{Z}$, the series $A(t)^m$ has integer coefficients. 
 
A power structure is called finitely determined
if for every $N>0$ there exists such $M>0$ that
the $N$-jet of the series $(A(t))^{m}$ is uniquely determined by the
$M$-jet of the series $A(t)$.

To fix a finitely determined power structure
it is sufficient to define the series $(1-t)^{-m}$
for every $m\in R$ such that $(1-t)^{-m-n}=(1-t)^{-m}\cdot(1-t)^{-n}$.
On the Grothendieck ring of varieties there is 
the power structure, defined by the equation
$$(1-t)^{-[X]}=1+[S^{1}X]t+[S^{2}X]t^{2}+\ldots,$$
where $S^{k}X$ denotes the $k$-th symmetric power of the variety $X$.
For example, for $j\ge 0$ Lemma 1 implies that
$$(1-t)^{-\mathbb{L}^j}=\sum_{k=0}^{\infty}t^k\mathbb{L}^{kj}=(1-t\mathbb{L}^j)^{-1}.$$ Macdonald formula for homologies of symmetric powers implies that  
$$(1-t)^{-\chi(X)}=1+\chi(S^{1}X)t+\chi(S^{2}X)t^{2}+\ldots.$$
This means that the Euler characteristic defines a morphism of power structures:
$$\chi((A(t))^m)=\chi(A(t))^{\chi(m)}$$
for all $m\in K_0(Var_{\mathbb{C}})$ and $A(t)=1+tK_0(Var_{\mathbb{C}})[[t]]$.







\begin{proposition}(\cite{powers})
There exists a natural power structure over the ring $K_0(Var_{\mathbb{C}})[\mathbb{L}^{-1}]$.
It can be uniquely continued to the continuous power structure over
 $\mathcal{M}$.
\end{proposition}

\section{ Correspondence between functions and curves}

Let us construct a measure on the symmetric power $S^{k}\mathcal{L}$ of the space $\mathcal{L}$  of arcs on $(\mathbb{C}^2,0)$.

The symmetric power $S^{k}\mathcal{L}_{n}$ is a
quasiprojective variety and one can define a natural continuation of the projection
$\pi_{n}^{(k)}:S^{k}\mathcal{L}\rightarrow S^{k}\mathcal{L}_{n}.$ Subset $A\subset S^{k}\mathcal{L}$
will be said to be cylindric, if there exist $n$ and a constructible subset $A_n\subset S^{k}\mathcal{L}_{n}$ such that
$A=(\pi_n^{(k)})^{-1}(A_n)$ . For this case let $\mu(A):=[A_n]\cdot\mathbb{L}^{-2kn}.$ Since  $[S^{k}\mathcal{L}_{n}]=\mathbb{L}^{-2kn}$, this definition does not depend of $n$.

For the cylindric set $D=\pi_{n}^{-1}(D_n)$ in $\mathcal{L}$
it is easy to see $\mu(S^{k}D)=\mathbb{L}^{-2nk}[S^{k}D_{n}]$.
Therefore the described construction corresponds to the power structure
over $\mathcal{M}$,
so that $$\sum_{k}\mu(S^{k}D)t^{k}=(1-t)^{-\mu(D)}.$$

Let $B_{i}$ be a collection of non-intersecting cylindric subsets of
$\mathcal{L}$, $k_i$ are nonnegative integers with 
 $\sum k_i=n$. Then for the natural embedding
$$S^{k_1}B_1\times S^{k_2}B_2\times\ldots\rightarrow S^{n}\mathcal{L}$$ 
it is easy to see that $\mu(\prod_{i}S^{k_i}B_{i})=\prod_{i}\mu(S_{k_i}B_{i}).$

It is clear that the measure of the set of tuples of arcs, where some of them coincide, is zero in the completed ring $\mathcal{M}$ since this set has infinite codimension.

\begin{lemma}(\cite{eqneng})
Let $f$ be a simple function on $\mathcal{L}$.
Define a function $F$ on $\sqcup_{k}S^{k}\mathcal{L}$
by the formula $F(\gamma_1,\ldots,\gamma_k)=\prod_{i}f(\gamma_{i}).$ 
Then
$$\int_{\sqcup_{k}S^{k}\mathcal{L}}Fd\chi_g=\int_{\mathcal{L}}(1-f)^{-d\chi_g}$$
if the right hand side integral converges.
Here $d\chi_g$ is in the exponent to emphasize that $1-f$ is considered as an element of an abelian group with respect to multiplication. 
\end{lemma}


In a similar way one can extend the motivic measure from the space
of branches $\mathcal{B}$ to the union $\sqcup_{k}S^{k}\mathcal{B}$
of its symmetric powers.
Let $p:\sqcup_{k}S^{k}(\mathcal{L}^{*}/\mathbb{C}^{*})\rightarrow\sqcup_{k}S^{k}\mathcal{B}$
be the natural projection.

\begin{lemma}
Let $F$ be a simple integrable function on 
$\sqcup_k S^{k}\mathcal{B}$. 
Then
$$\int_{\sqcup_{k}S^{k}(\mathcal{L}^{*}/\mathbb{C}^{*})}(p^{*}F)\cdot\mathbb{L}^{v}d\mu=
\int_{\sqcup_{k}S^{k}\mathcal{B}}F d\mu,$$
where $v$ denotes, as above, the total order of a tuple of arcs at the origin.
\end{lemma}
 


\begin{definition}
An arc $\gamma(t)=(x(t),y(t))$ is said to be degenerate, if there exist $h(x),$ $x^*(t),y^*(t)$ such that the order of $h$ is bigger or equal to 2,
and $x(t)=x^*(h(t)), y(t)=y^*(h(t)).$
\end{definition}

It is easy to prove that the set of degenerate curves has the motivic measure 0. This means, that up to the set of measure zero every arc can be considered as a uniformisation of its image.

Consider the natural map
$$Z:\mathbb{P}\mathcal{O}\rightarrow \sqcup_{k}S^k\mathcal{B},$$
which maps a function $f$ (defined up to multiplication by a constant) to the collection of uniformisations
of irreducible components of the germ of the curve 
$\{f=0\}$. If an irreducible component of $\{f=0\}$ has multiplicity
$m$, $Z(f)$ contains $m$ copies of its uniformisation.
It turns out that the transformation of the motivic measure under the map
$Z$ can be described explicitly.   

\begin{lemma} 
Let $f,g\in \mathcal{O}$, $g(\gamma(t))=0$ and $j_{n}(f)=j_n(g).$
Let $Ord_0(\gamma)=m,$ $\min\{Ord_0({\partial f\over\partial x}(\gamma(t))),Ord_0({\partial f\over\partial y}(\gamma(t)))\}=Q$ and $n>4Q$.
Consider  $n_1=mn-Q$.
There exists such arc $\gamma'(t)$,
that $j_{n_1}(\gamma')=j_{n_1}(\gamma)$ and
$f(\gamma'(t))=0.$
\end{lemma}

\begin{proof}
Note that $j_{mn}(f-g)(\gamma)=0, $
so $j_{mn}(f(\gamma))=0$.

Without any loss of generality, one can consider that
$ord(f_{y}(\gamma))=Q.$
This can be done using a generic affine change of variables
on the plane.

To construct an arc $\gamma'$ we will use the following process.
Let $$\gamma_{0}=\gamma,\gamma_{k+1}=\gamma_{k}-(0,
{f(\gamma_k)\over f_{y}(\gamma_k)}).$$

Then by Taylor formula 
$$f(\gamma_{k+1})=f(\gamma_{k}-{f(\gamma_k)\over f_{y}(\gamma_k)})=f(\gamma_k)-f_{y}(\gamma_k)\cdot {f(\gamma_k)\over f_{y}(\gamma_k)}+O(({f(\gamma_k)\over f_{y}(\gamma_k)})^2)=
O(({f(\gamma_k)\over f_{y}(\gamma_k)})^2).$$ Therefore if
$f(\gamma_k)$ is of order $N$ , then the order of
$f(\gamma_{k+1})$ is greater or equal to $2(N-Q)$. Since 
$f(\gamma_0)$
has the order greater than $mn$, the $(mn-Q)$-jets of
$\gamma_1$ and $\gamma_0$ coincide. Thus for all $k$
the $n_1$-jets of $\gamma_k$ and $\gamma$ are equal.
Therefore the order of $f(\gamma_1)$ is greater or equal to
$2(mn-Q)$, that is greater than ${3\over 2}mn$ if $n>4Q$.
We conclude that the sequence $\{\gamma_k\}$
converges in the $t$-adic topology, and the order of 
$f(\gamma_k)$ is greater or equal to $({3\over 2})^{k}mn$, so it
tends to infinity when $k$ tends to infinity.
Thus there exists a limit $\gamma'=\lim_{k\rightarrow\infty}\gamma_k$,
so $f(\gamma')=0$ and $j_{n_1}(\gamma')=j_{n_1}(\gamma)$.  
\end{proof}

\begin{lemma}
Let $\gamma_1,\ldots,\gamma_k$ be a collection of distinct arcs,
consider the linear map
$$ev_{\gamma}:\mathcal{O}\rightarrow \oplus_{i=0}^{k}\mathcal{O}_{\mathbb{C},0},
f\mapsto (f(\gamma_1(t)),\ldots,f(\gamma_k(t))).$$
Let $\mu$ be the Milnor number (see, e.g.,\cite{book}) of a function which determines the
union $\gamma$ of these arcs, $\delta(\gamma)={1\over 2}(\mu+k-1).$
Then the codimension of the image of $ev_{\gamma}$ is equal to $\delta(\gamma)$.
\end{lemma}

This lemma follows from the results of  \cite{codaira}.

For a (generally speaking, reducible) curve, defined by an equation
$f(x,y)=0$, with parametrisations $\gamma_1(t),\ldots,\gamma_k(t)$ of its components, consider
$$P_i(\gamma)=\min\{Ord_0({\partial f\over\partial x}(\gamma_i(t))),Ord_0({\partial f\over\partial y}(\gamma_i(t)))\}, P(\gamma)=\sum P_{i}(\gamma).$$

\begin{theorem}
Let $M\subset S^{k}\mathcal{B}$ be a measurable subset,
$N=Z^{-1}(M)$.
Then
$$\mu(N)=\int_{M}\mathbb{L}^{\delta(\gamma)-k-P(\gamma)}d\chi_{g},$$
where $\delta(\gamma)$ is defined in the previous lemma.
\end{theorem}
  
\begin{proof}
It is sufficient to suppose, that
 the orders of curves in $M$ are equal to $m_1,\ldots,m_k$,
then the order of every function  $f\in N$ is equal to $m=\sum_{i=1}^{k}m_i$. 
Without loss of generality one can suppose that
 $\delta(\gamma)$ is constant on  $M$ as well.

Consider the map $\zeta$ from the space of $n$-jets of functions
to the space of tuples of $(m_{i}n-P_i)$-jets
of arcs, transforming the polynomial $f$
to the collection of the corresponding jets of  $Z(f)$. 
Let us prove that for $n$ sufficiently large
its image coincides with the set $\pi(M)$ 
of corresponding collections of jets of curves from $M$.

If $(\gamma_1,\ldots,\gamma_k)=Z(f)$, $f_1=j_n(f)$,
then $n$-jets of $f$ and $f_1$ coincide, so
by Lemma 4 there exist such curves  $(\widetilde{\gamma}_1,\ldots,\widetilde{\gamma}_k),$
that $f_1(\widetilde{\gamma_i})=0$ and $j_{m_{i}n-P_i}(\widetilde{\gamma}_i)=j_{m_{i}n-P_i}(\gamma_i).$
Therefore $\zeta(f_1)$ coincides with the projection of $(\gamma_1,\ldots,\gamma_k)$.

If $f$ belongs to $N$, then $(\gamma_1,\ldots,\gamma_k)$
belongs to $M$, thus $\zeta(\pi_n(f))$ belongs to
$\pi(M)$. Inversely, if $(\gamma_1,\ldots,\gamma_k)\in M$,
and $\{f=0\}$ determines this collection of curves,
then $\pi(\gamma_1,\ldots,\gamma_k)=\zeta(\pi_{n}(f))\in\zeta(\pi_n(N)).$
So, $$\zeta(\pi_{n}(N))=\pi(M).$$

Let us describe fibers of the map $\zeta$. From the proof of Lemma 4 it is clear that  $j_{mn}(f(\gamma_0))=0$ if and only if there exists such arc $\widetilde{\gamma}_0$ that  $f(\widetilde{\gamma}_0)=0$
and $j_{n_1}(\gamma_0)=j_{n_1}(\widetilde{\gamma}_0).$

Let $\gamma=(\gamma_1,\ldots,\gamma_k)\in \pi(M).$
Let us calculate the dimension of the space $E_n$ of such $n$-jets $f$,
that for every $i$ $j_{m_{i}n}(f(\gamma_i))=0.$
Consider the linear map 
$$ev_{\gamma}:f\mapsto (j_{m_{i}n}(f(\gamma_i)).$$
By Lemma 5 the codimension of the image of this map is equal to
$\delta(\gamma)-k+1$ for $n$ sufficiently large (since we consider only series without degree 0 terms), 
so the codimension of its kernel is equal to the dimension of the 
image $\sum_{i}m_{i}n-\delta(\gamma)+k-1=mn-\delta(\gamma)+k-1.$    
In the space $E_n$ consider the subspace $E'_n$
of functions with the order greater than $m$. 
Since all functions with order $m$ 
vanishing on $\gamma$ have the same principal part
(homogeneous summand with degree $m$) up to multiplication by a constant, then the quotient $E_{n}/E'_n$ is one-dimensional,
so $\dim E'_n={(n+1)(n+2)\over 2}-1-(mn-\delta(\gamma)+k-1)-1.$ On the other hand,
if $f\in \zeta^{-1}(\gamma),$ then the fiber mapping to $\gamma$ is equal to 
$$\zeta^{-1}(\gamma)=\{\lambda f+h|\lambda\neq 0, h\in E'_n\},$$
so $[\zeta^{-1}(\gamma)]=\times[E'_n]$ (factor $\mathbb{L}-1$ corresponding to the parameter $\lambda$ is omitted because the space of functions is projectivized).
Thus in the Grothendieck ring one has
$$[\pi_{n}(N)]=\mathbb{L}^{{(n+1)(n+2)\over 2}-(mn-\delta(\gamma)+k-1)-2}[\pi(M)]$$
 
and therefore
$$\mathbb{L}\cdot\mathbb{L}^{-{(n+1)(n+2)\over 2}}[\pi_{n}(N)]=
\mathbb{L}\cdot\mathbb{L}^{-{(n+1)(n+2)\over 2}}\cdot\mathbb{L}^{{(n+1)(n+2)\over 2}-(mn-\delta(\gamma)+k-1)-2}[\pi(M)]=$$
$$\mathbb{L}^{\delta(\gamma)-k-mn}[\pi(M)]=
\mathbb{L}^{delta(\gamma)-k-P}\mu(M).$$
The last equation follows from 
$$\mu(M)=[\pi(M)]\cdot
\mathbb{L}^{-\sum_{i}(m_{i}n-P_i)}=[\pi(M)]\cdot
\mathbb{L}^{-mn+P}.$$
\end{proof}

Thus on the space of branches one can introduce the universal
factor
$$R(\gamma)=\mathbb{L}^{\delta(\gamma)-k-P(\gamma)}$$ such that $$\mu(Z^{-1}(M))=\int_{M}R(\gamma)d\mu,$$
so for every simple function $a(\gamma)$ 
$$\int_{\mathbb{P}\mathcal{O}}a(Z(f))d\mu=\int_{\sqcup_{k}S^{k}\mathcal{B}}R(\gamma)a(\gamma)d\mu.$$
Note that this equation yields that
$$\int_{\mathcal{O}}a(Z(f))d\mu=(\mathbb{L}-1)\int_{\sqcup_{k}S^{k}\mathcal{B}}R(\gamma)a(\gamma)d\mu.$$

{\bf Example 1.} Let $M_{k}\subset S^{k}\mathcal{B}$ be the set of $k$-tuples of smooth arcs with pairwise different tangents at the origin, $N_k=Z^{-1}(M_k)$. Let $(S^{k}\mathbb{P}^{1})^{*}$ be the set of unordered $k$-tuples of distinct points on $\mathbb{P}^{1}$.
It is clear that 
$$\mu(M_{k})=[(S^{k}\mathbb{P}^{1})^{*}]\cdot\mathbb{L}^{-k}.$$
The corresponding set $N_k\subset \mathcal{O}$ contains functions of order $k$ such that their homogeneous summand of order $k$ has $k$ different roots.
So one has
$$\mu(N_{k})=(\mathbb{L}-1)\cdot[(S^{k}\mathbb{P}^{1})^{*}] \cdot\mathbb{L}^{1-{(k+1)(k+2)\over 2}}.$$
In this case $\delta(\gamma)={k(k-1)\over 2}, P_i=k-1, P(\gamma)=k(k-1)$ so the proposition of Theorem 1
can be checked explicitly.   
 
{\bf Example 2.} Let $M_{2k}\subset \mathcal{B}$ -- be the set of arcs with the singularity of type $A_{2k}$ at the origin.
Suppose that an arc is tangent to the $x$-axis. Then 
one can choose a parametrisation such that $x(t)=t^2, y(t)=\sum_{k=1}^{\infty} y_k t^k$.
Then an arc has a singularity of type  $A_{2k}$  if and only if $y_2=0, y_1=y_3=\ldots=y_{2k-1}=0, y_{2k+1}\neq 0$, therefore the measure of the set of such arcs equals  $(\mathbb{L}-1)\mathbb{L}^{-k-2}$.
Since the tangent could be an arbitrary element of  $\mathbb{P}^1$, we have
 $$\mu(M_{2k})=(\mathbb{L}+1)(\mathbb{L}-1)\mathbb{L}^{-k-2}.$$
In this case $\delta(\gamma)=k, P(\gamma)=2k+1$. Therefore, 
$$\mu(N_{2k})=(\mathbb{L}+1)(\mathbb{L}-1)^2\mathbb{L}^{-2k-4}.$$

Let $k>1, M_{2k-1}\subset S^2\mathcal{B}$ be the set of arcs with the singularity of type $A_{2k-1}$ at the origin.
Suppose that an arc is tangent to the $x$-axis. Then 
one can choose a parametrisation such that $x^{(1)}(t)=t, x^{(2)}(t)=t$.
The singularity is of type $A_{2k-1}$ if and only if $Ord_{0}(y^{(1)}(t)-y^{(2)}(t))=k$, so a measure of the set of such arcs equals to $(\mathbb{L}-1)\mathbb{L}^{-k}$.
Since tangent varies over the projective line $\mathbb{P}^1$, one has
 $$\mu(M_{2k-1})=(\mathbb{L}+1)(\mathbb{L}-1)\mathbb{L}^{-k-1}.$$
In this case $\delta(\gamma)=k, P(\gamma)=2k$. Thus, 
$$\mu(N_{2k-1})=(\mathbb{L}+1)(\mathbb{L}-1)^2\mathbb{L}^{-2k-3}.$$

These answers can be checked explicitly if one considers the 
requirement on a function to have a prescribed singularity $A_s$
in terms of its Newton diagram.

From the previous example we know that the measure of the set of functions with the singularity $A_1$ equals to $$\mu(A_1)=(\mathbb{L}-1)\mathbb{L}^{-3}.$$

Let us calculate the sum of the measures of the sets of functions of type $A_m$ over all  $m$:
$$\mu(A)=(\mathbb{L}-1)\mathbb{L}^{-3}+\sum_{k=1}^{\infty}(\mathbb{L}+1)(\mathbb{L}-1)^2\mathbb{L}^{-2k-4}+\sum_{k=2}^{\infty}(\mathbb{L}+1)(\mathbb{L}-1)^2\mathbb{L}^{-2k-3}$$ $$=(\mathbb{L}-1)\mathbb{L}^{-3}+(\mathbb{L}+1)(\mathbb{L}-1)^2 ({\mathbb{L}^{-6}\over 1-\mathbb{L}^{-2}}+{\mathbb{L}^{-7}\over 1-\mathbb{L}^{-2}})=(\mathbb{L}-1)(\mathbb{L}^{-3}+\mathbb{L}^{-4}+\mathbb{L}^{-5})=\mathbb{L}^{-2}-\mathbb{L}^{-5}.$$

On the other hand,  the function has $A_m$ singularity  if and only if its 1-jet is zero, and 2-jet is nonzero, so (considering the space of 2-jets):
$$\mu(A)=\mathbb{L}^{-5}(\mathbb{L}^{3}-1),$$
what is equal to the previous answer.

{\bf Example 3.} Let us calculate the measure of the $\{\mu=const\}$ stratum in $\mathcal{O}$ for the singularity  $x^p+y^q$ for coprime $p$ and $q$.

From the parametrisation viewpoint, if an arc is tangent to the $x$-axis, one can consider  $x(t)=t^{p}$. An arc is in the prescribed stratum, if the smallest number of monomial in $y(t)$, which is not multiple of  $p$  and its coefficient   does not vanish, is $q$, and coefficient at $t^p$ vanishes. The measure of the set of such arcs is equal to $(\mathbb{L}-1)\mathbb{L}^{-q+[{q\over p}]-1},$
therefore $$\mu(M)=(\mathbb{L}+1)(\mathbb{L}-1)\mathbb{L}^{-q+[{q\over p}]-1}.$$
One has $\delta(\gamma)={(p-1)(q-1)\over 2}, P(\gamma)=(p-1)q,$
thus 
$$\mu(N)=(\mathbb{L}+1)(\mathbb{L}-1)\mathbb{L}^{[{q\over p}]-1-{(p+1)(q+1)\over 2}}.$$
For example, the codimension of the set $N$ is equal to 
$$c={(p+1)(q+1)\over 2}-2-[{q\over p}].$$
One has a formula (\cite{book}), combining the codimension of the $\{\mu=const\}$ stratum, the Milnor number of a singularity and the modality $m$:
$$\mu=c+m-1,$$
so from the known codimension we can obtain the following formula for the modality:
$$m={pq\over 2}-{3p\over 2}-{3q\over 2}+3,5+[{q\over p}].$$
One can check that this expression coincides with the Kouchnirenko's formula (\cite{kushn}), describing modality of this singularity as the number of integer points in the rectangle with vertices   $(0,0),(0,p-2),(q-2,0),(q-2,p-2)$, situated above the line $px+qy=pq$.

The invariant $P(\gamma)$ of a curve singularity can be expressed via its Milnor number.

\begin{lemma}
Let $\mu$ and $v$ be the Milnor number and the vanishing order of a curve  $\gamma$. Then
\begin{equation}
\label{pmu}
P(\gamma)=\mu+v-1.
\end{equation} 
\end{lemma}
\begin{proof}
Suppose at first that the curve $\gamma$ is irreducible with a parametrisation
$(x(t),y(t))$ and an equation $f=0$. Without  loss of generality one can assume that the order of $x(t)$ is less than the one of $y(t)$.  Since
$f(x(t),y(t))\equiv 0,$  $$f_{x}(\gamma(t))\dot{x}(t)+f_{y}(\gamma(t))\dot{y}(t)=0.$$
Vanishing orders of summands are equal. Therefore
$$Ord\, f_{x}(\gamma(t))+v_x-1=Ord\, f_{y}(\gamma(t))+v_y-1,$$
thus $$Ord\, f_x(\gamma(t))=Ord\, f_{y}(\gamma(t))+(v_y-v_x),$$
and $P(\gamma)=Ord\, f_{y}(\gamma(t))$, 
$$P(\gamma)-v=Ord\, f_{y}(\gamma(t))-v_x=Ord\, f_{x}(\gamma(t))-v_y.$$

Consider the blow-up of the origin. An equation of the strict transform of the curve is
$\hat{f}=x^{-v}f(x,\theta x)=0$, and 
$$\hat{P}-\hat{v}=Ord\,\hat{f}_{\theta}-v_x=Ord\, [x^{1-v}f_{y}(\gamma(t))]-v_x=P-v-v(v-1).$$
On the other hand, it is easy to check that
$\hat{\mu}=\mu-v(v-1),$  so if the proposition is true for the strict transform, it is also true for the initial curve. Therefore it should be checked only for nonsingular curve, where  $P=0, v=1, \mu=0,$ so the proposition of lemma is true.

Let us prove now the statement of the lemma by induction on the number of components of a curve.
Let $\gamma=\gamma_1\cup\gamma_2$, then $f=f_1f_2$.
For every component of $\gamma_1$, the order of $f_x$ on it equals to the sum of the orders of $f_{1x}$ and $f_2$, therefore the sum of such orders over the components of $\gamma_1$ equals to $P(\gamma_1)+\gamma_1\circ\gamma_2$. The same statement is true for $\gamma_2$ and after addition we have
$$P(\gamma)=P(\gamma_1)+P(\gamma_2)+2\gamma_1\circ\gamma_2.$$
On the other hand, 
one has
$$\mu(\gamma)=\mu(\gamma_1)+\mu(\gamma_2)+2\gamma_1\circ\gamma_2-1.$$
By the induction assumption (\ref{pmu}) is true for $\gamma_1$ and $\gamma_2$, therefore it is true for $\gamma$ also. 
\end{proof}

\begin{theorem}
Let $a$ be a measurable function on $\sqcup_{k}S^{k}\mathcal{B}.$ Then
$$\int_{\mathbb{P}\mathcal{O}}a(Z(f))d\mu=
\int_{\sqcup_{k}S^{k}(\mathcal{L}^{*}/\mathbb{C}^{*})}
(p^{*}a)\mathbb{L}^{-\delta(\gamma)}d\mu.$$
\end{theorem}

\begin{proof}
By Theorem 1,
$$\int_{\mathcal{O}}a(Z(f))d\mu=
\int_{\sqcup_{k}S^{k}\mathcal{B}}
a\mathbb{L}^{\delta(\gamma)-k-P(\gamma)}d\mu=
\int_{\sqcup_{k}S^{k}(\mathcal{L}^{*}/\mathbb{C}^{*})}
(p^{*}a)\mathbb{L}^{\delta(\gamma)-k-\mu-v+1}\mathbb{L}^{v}d\mu.$$
Now we can use the equation
$\delta(\gamma)={1\over 2}(\mu+k-1).$
\end{proof}

It is important to note that $\delta(\gamma)$ has a simple geometric meaning:
it is the number of self-intersections of a generic deformation of $\gamma$.

{\bf Example 4.} $$(\mathbb{L}-1)\int_{\sqcup_{k}S^{k}(\mathcal{L}^{*}/\mathbb{C}^{*})}a^{v_x}b^{v_y}\mathbb{L}^{-\delta(\gamma)}d\mu=
\int_{\mathcal{O}}a^{v_x}b^{v_y}d\mu={ab\mathbb{L}^{-2}(\mathbb{L}-1)^2\over (1-a\mathbb{L}^{-1})(1-b\mathbb{L}^{-1})}.$$
The author do not know the simple explanation of this fact.
For example, the same integral over the space $\mathcal{L}^{*}/\mathbb{C}^{*}$ itself is much more complicated and satisfy some functional equations (\cite{eqneng}). One can say that addition of all other symmetric powers "simplify" the integral.

{\bf Example 5.}
Let $\gamma=(\gamma_1,\ldots,\gamma_k)$ be a tuple of arcs on the plane at the origin. It provides a $\mathbb{Z}^k$-valued function on the space $\mathcal{O}$:
$$\underline{v}_{\gamma}(f)=(Ord_{0}f(\gamma_1(t)),\ldots,Ord_{0}f(\gamma_k(t))).$$

Consider the integral
$$P(\underline{t})=\int_{\mathbb{P}\mathcal{O}}\underline{t}^{\underline{v}_{\gamma}}\mathbb{L}^{\delta(f)}d\mu.$$

If $\omega$ is a collection of uniformisations of the components of  $\{f=0\}$, then $\underline{v}_{\gamma}(f)$ equals to the tuple of intersection numbers of $\gamma_k$ and $\omega$.
So Theorem 2  yields 
$$P(\underline{t})=\int_{\sqcup_{k}S^{k}(\mathcal{L}^{*}/\mathbb{C}^{*})}
\underline{t}^{\underline{\gamma}\circ\gamma_1}d\mu=
\int_{\mathcal{L}^{*}/\mathbb{C}^{*}}(1-\underline{t}^{\underline{\gamma}\circ\gamma_1})
^{-d\mu}.$$
Last equality follows from Lemma 2.

Consider an embedded resolution of an arc $\gamma$.
It is a proper map $(X,D)\rightarrow (\mathbb{C}^2,0)$, which is an isomorphism out of origin, such that 
the preimage of $\gamma$  is a normal crossing divisor.
The exceptional divisor is a normal crossing divisor $D$, which
components $E_s$ are isomorphic to $\mathbb{CP}^1$. 
Let $J=\sum_{s}\nu_{s}E_s$ be the relative canonical divisor
(locally given by the Jacobian).
Let $(m_{ij})$ be the inverse matrix to $(-(E_{i}\circ E_{j}))$.

Let $$F(q)=\prod_{k}{1\over 1-q^k},\qquad G(x,y)=\prod_{k,m}{1\over 1-x^ky^m} $$

The arguments of (\cite{nonpar}, Theorem 1) adapted to the 2-dimensional case 
yields the following equation:

$$P(\underline{t})=\prod_{s}F(\underline{t}^{\underline{m}_s}\mathbb{L}^{-(\nu_s+1)})^{-[E_s^{\circ}]}
\prod_{E_{i}\cap E_{j}\neq\emptyset}G(\underline{t}^{\underline{m}_i}\mathbb{L}^{-(\nu_i+1)},\underline{t}^{\underline{m}_j}\mathbb{L}^{-(\nu_j+1)})^{-(\mathbb{L}-1)}\times$$
$$\prod_{E_{i}\cap\widetilde{\gamma}_{j}\neq \emptyset}G(\underline{t}^{\underline{m}_i}\mathbb{L}^{-(\nu_i+1)},\underline{t}^{\underline{1}_j}\mathbb{L}^{-1})^{-(\mathbb{L}-1)}.
$$

\section{Correspondence for the Euler characteristic}

Since the set of motivic measure zero can have non-zero Euler characteristic, one cannot use naively the result of the Theorem 2 applied to the Euler characteristic. But in some cases there are nice formulas for such correspondence.

Let $\mathcal{U}\subset \mathcal{B}$ be the set of nondegenerate arcs, i. e. arcs which are uniformisations of their images.

It is clear that the image of the map $Z:\mathbb{P}\mathcal{O}\rightarrow \sqcup_{k}S^{k}\mathcal{B} $ coincides with $\sqcup_{k}S^{k}\mathcal{U}$ and fibers of this map are isomorphic to affine spaces.

Suppose that $a(\gamma)$ is such  a function on  $\mathcal{L}$ that if $\gamma(t)=\gamma'(h(t))$ then $$a(\gamma)=(a(\gamma'))^{Ord_{0}h(t)}.$$  For example, it is an $Aut_{\mathbb{C},0}$ -- invariant. Let $A$ be the function on $\sqcup_{k}S^{k}\mathcal{L}$ defined by the formula
$$A(\gamma_1,\ldots,\gamma_k)=a(\gamma_1)\cdot\ldots\cdot a(\gamma_k).$$

\begin{theorem}
$$\int_{\mathcal{L}^{*}/\mathbb{C}^{*}}(1-a(\gamma))^{-d\chi}=
\prod_{k=1}^{\infty}\int_{\mathbb{P}\mathcal{O}}A(Z(f))^{k}d\chi.$$
\end{theorem}
\begin{proof}
Since the Euler characteristics of fibers of the map $Z$ are equal
to 1, by the Fubini formula one has $$\prod_{k=1}^{\infty}\int_{\mathbb{P}\mathcal{O}}A(Z(f))^{k}d\chi=\prod_{k=1}^{\infty}\int_{\sqcup_{m}S^{m}\mathcal{U}}A(\gamma)^{k}d\chi=\prod_{k=1}^{\infty}\int_{\mathcal{U}}(1-a(\gamma)^{k})^{-d\chi}.$$ On the other hand, every  arc $\gamma$ from $\mathcal{L}^{*}/\mathbb{C}^{*}$ has a unique 
uniformisation $\gamma'$ up to automorphism of $(\mathbb{C},0)$, so $\gamma=\gamma'(h(t))$. Since we factorize by $\mathbb{C}^{*}$, one can assume that
 $h(t)=t^{k}+\ldots.$ Hence $(\mathcal{L}^{*}/\mathbb{C}^{*})$ is decomposed to parts with fixed order of $h(t)$, and the projection of each part onto $\mathcal{U}$ has affine fibers. Therefore  
$$\int_{\mathcal{L}^{*}/\mathbb{C}^{*}}(1-a(\gamma))^{-d\chi}=
\prod_{k=1}^{\infty}\int_{\mathcal{U}}(1-a(\gamma)^{k})^{-d\chi}.$$
\end{proof}
\begin{corollary}
Let $\mu$ be the Moebius function. Then $$\int_{\mathbb{P}\mathcal{O}}A(Z(f))d\chi=\prod_{k=1}^{\infty}
\int_{\mathcal{L}^{*}/\mathbb{C}^{*}}(1-a^{k})^{-\mu(k)d\chi}.$$
\end{corollary}
\begin{proof}
By Theorem 3
$$\prod_{k=1}^{\infty}\int_{\mathcal{L}^{*}/\mathbb{C}^{*}}(1-a^{k})^{-\mu(k)d\chi}=\prod_{k,m=1}^{\infty}(\int_{\mathbb{P}\mathcal{O}}A^{km}(Z(f))d\chi)^{\mu(k)},$$
so the proposition follows from the Moebius inversion formula .
\end{proof}

\section*{Acknowledgments}

I am grateful to  my advisor S. Gusein-Zade, who pushed me to search for such correspondence,  for constant attention and
encouragement.

Moscow State University,\newline
E.mail: gorsky@mccme.ru


\begin{thebibliography}{99}

\bibitem{book}
V. I. Arnold, S. M. Gusein-Zade, A. N. Varchenko. 
Singularities of differentiable maps. Vol. 2,
Birkhauser, 1985.

\bibitem{al}
A. Campillo, F. Delgado, S. M. Gusein-Zade.
The Alexander polynomial of a plane curve singularity via the ring of functions on it.
Duke Math J. 117 (2003), no. 1, 125--156.

\bibitem{al2}
A. Campillo, F. Delgado, S. M. Gusein-Zade.
Integrals with respect to the Euler characteristic over spaces of functions and the Alexander polynomial.
Proc. Steklov Inst. Math. 2002, no. 3 (238), 134--147.

\bibitem{alex}
A. Campillo, F. Delgado, S. M. Gusein-Zade. 
Multi-index filtrations and motivic Poincar\'e series.
arXiv: math.AG/0406240 (to appear in Monatshefte f$\ddot{u}$r Mathematik).


\bibitem{dl}
J. Denef, F. Loeser. Germs of arcs on singular algebraic varieties and 
motivic integration. Inventiones Math. 135 (1999), no.1, 201-232. 

\bibitem{eqneng}
E. Gorsky.
Motivic integrals and functional equations. 
arXiv: math.AG/0606521 


\bibitem{powers}
S. M. Gusein-Zade, I. Luengo, A. Melle-Hern\'andez.
A power structure over the Grothendieck ring of varieties.
Math. Res. Lett. 11(2004), no.1, 49-57. 

\bibitem{nonpar}
S. M. Gusein-Zade, I. Luengo, A. Melle-Hern\'andez.
Integration over spaces of nonparametrized arcs and motivic 
versions of the monodromy zeta function. 
Proceedings of the Steklov Institute of Mathematics, 2006, Vol. 252, pp. 63-73. 




\bibitem{codaira}
K. Kodaira. On compact complex analytic surfaces.
$I$. Ann. of Math. (2) 71, 1960, 111--152.

\bibitem{kushn}
A. G. Kouchnirenko. Polyederes de Newton et nombres de Milnor. - Invent. Math., 1976, v. 32, pp. 1-31.


\end{thebibliography}
\end{document}